\documentclass[12pt, twoside, reqno]{article}
\usepackage{amsmath,amsthm}
\usepackage{amssymb}
\usepackage{enumitem}

\newtheorem{theorem}{Theorem}[section]

\newtheorem{lemma}[theorem]{Lemma}

\theoremstyle{definition}

\newtheorem{remark}[theorem]{Remark}

\newtheorem{conjecture}[theorem]{Conjecture}

\numberwithin{equation}{section}

\begin{document}

\baselineskip=17pt


\title{Primes of the form $ax+by$}

\author{Yong-Gao Chen\footnote{Corresponding author} and Hui Zhu\\
\small School of Mathematical Sciences, Nanjing Normal University \\
\small  Nanjing  210023,  P.R. China\\
\small ygchen@njnu.edu.cn, zh200109150528@163.com}
\date{}

\maketitle

\begin{abstract} For two coprime positive integers $a,b$, let
$T(a,b)=\{ ax+by : x,y\in \mathbb{Z}_{\ge 0} \} $ and let
$s(a,b)=ab-a-b$. It is well known that all integers which are
greater than $s(a,b)$ are in $T(a,b)$. Let $\pi (a, b)$ be the
number of primes  in $T(a,b)$ which are less than or equal to
$s(a,b)$. It is easy to see that $\pi (2, 3)=0$ and $\pi (2, b)=1$
for all odd integers $b\ge 5$. In this paper, we prove that if
$b>a\ge 3$ with $\gcd (a, b)=1$, then $\pi (a, b)>0.005
s(a,b)/\log s(a,b)$. We conjecture that $\frac{13}{66}\pi
(s(a,b))\le \pi (a, b)\le \frac 12\pi (s(a,b))$ for all $b>a\ge 3$
with $\gcd (a, b)=1$.
\end{abstract}

\renewcommand{\thefootnote}{}



\renewcommand{\thefootnote}{\arabic{footnote}}
\setcounter{footnote}{0}

\section{Introduction}

For two coprime positive integers $a,b$, let $s(a,b)=ab-a-b$ and
$T(a,b)=\{ ax+by : x,y\in \mathbb{Z}_{\ge 0} \} $. It is well
known that if $a,b$ are two coprime positive integers, then all
integers which are greater than $s(a,b)$ are in $T(a,b)$ and
$s(a,b)\notin T(a,b)$. Let $\pi (a, b)$ be the number of primes in
$T(a,b)$ which are less than $s(a,b)$ and  $\pi (x)$ the number of
primes not exceeding $x$. Ram\'irez Alfons\'in and Ska\l ba
\cite{Ramiez2020} proved that for any $\varepsilon
>0$, there is a positive constant $C(\varepsilon)$ such that for $ab$ sufficiently
large, \begin{equation}\label{eq1}\pi (a,b)>C(\varepsilon)
\frac{s(a,b)}{(\log s(a,b))^{2+\varepsilon}}.\end{equation}
Ram\'irez Alfons\'in and Ska\l ba \cite{Ramiez2020} also posed the
following conjectures:

{\bf Conjecture A.} {\it Let $b>a\ge 3$ with $\gcd (a,b)=1$. Then
$$\pi (a, b)\sim \frac 12 \pi (s(a,b)) \text{ as } a\to \infty
.$$}

{\bf Conjecture B.} {\it Let $b>a\ge 3$ with $\gcd (a,b)=1$. Then
$$\pi (a, b)>0.$$}

Ding \cite{Ding2023} showed that Conjecture A is true for almost
all $a, b$ with $\gcd (a, b)=1$. Recently, Ding, Zhai and Zhao
\cite{DingZhaiZhaoarXiv2023} confirmed Conjecture A and Dai, Ding
and Wang \cite{DaiDingWangarXiv2024} confirmed Conjecture B.

In this paper,  the following results are proved. Our proofs do
not depend on the results in \cite{DaiDingWangarXiv2024}
 and \cite{DingZhaiZhaoarXiv2023}.

\begin{theorem}\label{thm1} Let $a\ge 3$ be a fixed integer. Then as $b\to \infty$ with $\gcd (a, b)=1$, we have
$$\pi (a, b)=\left( \frac 12-\frac 1{2(a-1)} +o(1)\right) \pi (s(a,b)).$$
In particular, $$\pi (3, b)=\frac 14 \pi (s(3,b))+o(\pi (s(3,b))
\text{ as } b\to \infty , \ 3\nmid b.$$
\end{theorem}

\begin{theorem}\label{thm2} If  $3\le a<b$ with $\gcd (a,b)=1$, then \begin{equation*}\pi
(a,b)>0.005 \frac {s(a,b)}{\log s(a,b)} .\end{equation*}
\end{theorem}

\begin{remark} Since $s(2, b)=b-2$, it follows that if $x,y\in \mathbb{Z}_{\ge 0}$ and $2x+by<s(2, b)$,
then $y=0$. So $\pi (2, b)=1$ for all odd integers $b\ge 5$. Thus,
we always assume that $a\ge 3$. It is clear that
\begin{equation*}\pi (a,b)\le \pi (s(a,b))=(1+o(1)) \frac {s(a,b)}{\log
s(a,b)}.\end{equation*} From the proof of Theorem \ref{thm2}, we
see that $0.005$ can be improved.\end{remark}

Basing on numerical calculations, we pose the following
conjectures.

\begin{conjecture} For $3\le a<b$ with $\gcd (a,b)=1$, we have
$$\pi (a, b) \ge \frac{13}{66} \pi (s(a,b))$$
and the equality holds if and only if $a=3$ and $b=166$.
\end{conjecture}

\begin{conjecture} For $3\le a<b$ with $\gcd (a,b)=1$, we have
$$\pi (a, b) \le \frac{1}{2} \pi (s(a,b))$$
and the equality holds if and only if $a=3$ and $b=5$.
\end{conjecture}

In the following, let $\mathcal{P}$ be the set of all primes,
$S=s(a,b)$ and $c=\frac 18 S$.

\section{Preliminary lemmas}

Let $\pi (x; m, l)$ denote the number of primes $p\le x$ with
$p\equiv l\pmod{m}$ and let $\varphi (n)$ be the Euler totient
function. For two coprime integers $b>a\ge 3$ and any nonnegative
integer $k$, let
$$S_k (a, b)= \{ xa+kb : x\in \mathbb{Z}_{\ge 0}, 0<xa+kb\le s(a,b)\}.$$

\begin{lemma}\cite[Theorem 2]{MontgomeryVaughan1973} \label{lem7}
Let $k,l$ be two integers and let $x,y$ be positive real numbers
with $1\le k<y$. Then
$$\pi (x+y; k,l)-\pi (x; k,l)<\frac{2y}{\varphi (k) \log (y/k)}.$$
\end{lemma}

\begin{remark} In \cite[Theorem 2]{MontgomeryVaughan1973},
it is required that $l$ is a positive integer. It is clear that
$l$ can be an arbitrary integer which is coprime to $k$.
\end{remark}

\begin{lemma}\cite[(1.15)]{Illinois2018}\label{PiEstimate} For $x\ge 1474279333$, we have
$$|\pi (x)-\mbox{Li} (x)|<0.0008375\frac{x}{\log^2 x}, $$
where
$$\mbox{Li} (x)=\int_2^x \frac 1{\log t} dt.$$
\end{lemma}

\begin{lemma}\label{PiEstimate8numberic} For $x\ge 40000$, we have
\begin{equation*}\pi (x)-\pi (\frac 78 x)>0.121243 \frac {x}{\log x}.\end{equation*}
\end{lemma}

\begin{proof} For $x\ge 10^{10}$, we have $\frac 78
x>1474279333$. By Lemma \ref{PiEstimate}, we have
\begin{align*}&\pi (x)-\pi (\frac 78 x)\\
&>\mbox{Li}(x)-0.0008375\frac{x}{\log^2 x} -\mbox{Li}(\frac
78x)-0.0008375\frac{7x/8}{\log^2 (7x/8)}\\
&=\int_{7x/8}^x \frac 1{\log t} dt-0.0008375\frac{x}{\log^2
x}\\
&\ \ \ -0.0008375\cdot \frac 78\cdot \frac 1{(1+(\log (7/8))/\log
x)^2}
\frac{x}{\log^2 x}\\
&>\int_{7x/8}^x \frac 1{\log t} dt-0.001579613\frac{x}{\log^2 x}\\
&>\frac 18 \frac {x}{\log x}-0.0015789\frac{x}{\log^2 x}\\
&>0.12493 \frac {x}{\log x}.
\end{align*}
For $10^{5}\le x< 10^{10}$, there exist two integers $u,k$ with
$2\le u\le 6$ and  $10^3\le k<10^4$ such that $10^uk\le
x<10^u(k+1)$. By a  calculation, we have \begin{align*}&\ \ \ \ (\pi (x)-\pi (\frac 78 x))\frac {\log x}{x}\\
&\ge (\pi (10^uk)-\pi (\frac 78 \cdot 10^u(k+1) ))\frac {\log (10^u(k+1))}{10^u(k+1)}\\
&> 0.122583.
\end{align*}
For $40000\le x< 10^{5}$, there exists an integer $40000\le k<
10^{5}$ with $k\le x<k+1$. By a  calculation, we have
$$(\pi (x)-\pi (\frac 78 x))\frac {\log x}{x}\ge (\pi (k)-\pi (\frac 78 (k+1)))\frac {\log (k+1)}{k+1} >0.121243.$$

This completes the proof of Lemma
\ref{PiEstimate8numberic}.\end{proof}

\begin{lemma} (Siegel--Walfisz Theorem, see \cite[Chapter 22]{Davenport1980}) \label{pi(x;m,l)=} Let $A$ be a positive real number and $m,l$ two coprime
integers with $1\le m\le (\log x)^A$. Then there is a positive
constant $D$ depending only on $A$ such that
$$\pi (x;m,l)=\frac 1{\varphi (m)} \int_2^x \frac 1{\log t} dt +O(x \exp (-D \sqrt{\log x})), $$
uniformly in $m$.
\end{lemma}

\begin{lemma}\cite[Corollary 1.6]{Illinois2018} \label{pi(x;m,l)estimatesmall} Let $m, l$ be two coprime integers with $1\le m\le 1200$. Then for all $x\ge 50m^2$,
$$\frac{x}{\varphi (m) \log x}< \pi (x;m,l)< \frac{x}{\varphi (m) \log x} \left( 1+\frac 5{2\log x}\right).$$
\end{lemma}

\begin{lemma}\label{lem1} Let $b>a\ge 3$ be two coprime integers.
Then $S_k (a, b)$ $(k=0, 1,\dots )$ are disjoint each other.
\end{lemma}

\begin{proof} Let $k_1<k_2$ be two  nonnegative integers.
Suppose that $S_{k_1} (a, b)\cap S_{k_2} (a, b)\not= \emptyset $.
Then there exist $x_1, x_2\in \mathbb{Z}_{\ge 0}$ such that
$$x_1a+k_1b=x_2a+k_2b\le  s(a,b).$$
So $(x_1-x_2)a=(k_2-k_1)b$. By $(a, b)=1$, we have $a\mid
k_2-k_1$. It follows that $k_2\ge k_1+a\ge a$. Thus, $x_2a+k_2b\ge
ab>s(a,b)$, a contradiction.

This completes the proof of Lemma \ref{lem1}.\end{proof}

\begin{lemma}\label{lem6} Let $m,k$ be two coprime positive integers.
Then for any positive real number $x$ and any integer $l$, we have
$$\sum_{1\le n\le x\atop (nk+l, m)=1} 1 =\frac{\varphi (m)}{m} x+\theta_{x, m, k, l} 2^{\omega (m)},$$
where $\omega (m)$ is the number of distinct prime divisors of $m$
and $\theta_{x, m, k, l}$ is a real number with $|\theta_{x, m, k,
l}|\le 1$.
\end{lemma}

\begin{proof} We have
\begin{align*}&\sum_{n\le x\atop (nk+l, m)=1} 1=\sum_{n\le x}
\sum_{d\mid nk+l\atop d\mid m} \mu (d)\\
&=\sum_{d\mid m} \mu (d) \sum_{n\le x\atop d\mid nk+l} 1\\
&=\sum_{d\mid m} \mu (d) \left( \frac xd +\delta_{x,d,k,l} \right)\\
&=\frac{\varphi (m)}{m} x+\theta_{x, m,k,l} 2^{\omega (m)},
\end{align*}
where $\delta_{x,d,k,l}$ and  $\theta_{x, m, k, l}$ are real
numbers with $|\delta_{x,d,k,l}|\le 1$ and $|\theta_{x, m, k,
l}|\le 1$.

This completes the proof of Lemma \ref{lem6}.\end{proof}

\begin{lemma}\label{lem9} Let $a,b$ be two integers with $b>a\ge 10$ and $\gcd (a, b)=1$. Then
\begin{align*}&\# \{ p\in \mathcal{P} : \frac 78S< p\le S, p\notin T( a, b)  \}\\
& < \frac 1{4\varphi (a)}\Big( \sum_{\substack{1\le n\le (a+6)/8
\\ \gcd (a, n)=1}} 1\Big) \Big( 1-\frac{\log (8a)}{\log
S}\Big)^{-1} \cdot \frac{S}{\log S}\end{align*} and
\begin{align*}&\# \{ p\in \mathcal{P} : \frac 78S< p\le S, p\notin T( a, b)  \} \\
&< \frac 1{4}\left( \frac18 + \frac{3}{4a}+ \frac{2^{\omega
(a)}}{\varphi (a)}\right) \left( 1-\frac{\log (8a)}{\log
S}\right)^{-1} \cdot \frac{S}{\log S}.\end{align*}
\end{lemma}

\begin{proof} By $b>a\ge 10$, we have $9a<11(a-1)\le b(a-1)$. So
$8a<b(a-1)-a=S$. That is, $a<\frac 18 S$. For $0\le p\le S$, we
have (see \cite{Ramiez2020})
\begin{align*}p\notin T( a, b) &\Longleftrightarrow S-p\in
T( a, b)  \Longleftrightarrow S-p=ax+by, x,y\in \mathbb{Z}_{\ge
0} \\
& \Longleftrightarrow p=S-ax-by, x,y\in \mathbb{Z}_{\ge
0}.\end{align*} Hence
\begin{align*}&\# \{ p\in \mathcal{P} : \frac 78S< p\le S, p\notin T( a, b)  \}\\
=& \# \{ p\in   \mathcal{P} : \frac 78S< p\le S, p= S-ax-by,
x,y\in
\mathbb{Z}_{\ge 0}, 0\le ax+by<\frac 18S \}\\
\le & \# \{ p\in   \mathcal{P} : \frac 78S< p\le S, p= S-ax-by,
x,y\in
\mathbb{Z}_{\ge 0}, 0\le by<\frac 18S \}\\
\le & \# \{ p\in  \mathcal{P} : \frac 78S< p\le S, p\equiv
-b-by\pmod{a}, y\in \mathbb{Z}, 0\le by< \frac S8 \} .
\end{align*}
Let $p\in  \mathcal{P}$ such that $\frac 78S< p\le S$, $p\equiv
-b-by\pmod{a}$ for some $y\in \mathbb{Z}$. Then $\gcd (a, 1+y)\mid
p$. If $\gcd (a, 1+y)>1$, then  $\gcd (a, 1+y)=p$ by $p\in
\mathcal{P}$. Thus, $p=\gcd (a, 1+y)\le a<\frac 18S$, a
contradiction with $p>\frac 78S$. Hence, $\gcd (a, 1+y)=1$.  If
$y$ is an integer with $0\le by< S/8$, then  $8by<
ab-a-b<b(a-1)$. So $8y\le a-2$. Thus, \begin{align*}&\# \{ p\in \mathcal{P} : \frac 78S< p\le S, p\notin T( a, b)  \} \\
&\le \sum_{\substack{0\le y\le (a-2)/8 \\ \gcd (a, 1+y)=1}} \# \{
p\in \mathcal{P} : \frac 78S< p\le S, p\equiv
-b-by\pmod{a} \}\nonumber\\
&= \sum_{\substack{1\le n\le (a+6)/8\\
\gcd (a,n)=1}} \left( \pi (S; a, -bn) -\pi (\frac 78S; a,
-bn)\right).
\end{align*}
By Lemma \ref{lem7} and $10\le a<\frac 18S$, we have
\begin{align*}
 &\pi (S; a,
-bn) -\pi (\frac 78S; a, -bn) \\ &<  \frac{S}{4\varphi
(a) \log (S/8a)}\\
&= \frac 1{4\varphi (a)}\left( 1-\frac{\log (8a)}{\log
S}\right)^{-1} \cdot \frac{S}{\log S}.
\end{align*}
By Lemma \ref{lem6}, we have
\begin{align*}\sum_{\substack{1\le n\le (a+6)/8\\
\gcd (a,n)=1}} 1 &\le \frac{\varphi (a)}{a}\frac{a+6}8 +
2^{\omega (a)}\\
&=\varphi (a) \left( \frac18 + \frac{3}{4a}+ \frac{2^{\omega
(a)}}{\varphi (a)}\right).
\end{align*}
Now Lemma \ref{lem9} follows immediately.\end{proof}

\begin{lemma}\label{lem8} Let $p_i$ be the $i$-th prime and $t$ a positive integer.
For any positive integer $n$ with $n\ge p_1\cdots p_t$, we have
\begin{align}\label{eq4}\frac{2^{\omega (n)}}{\varphi (n)}\le
\frac{2^t}{\varphi (p_1\cdots p_t)}.\end{align}

In particular, if $p_1\cdots p_t\le n<p_1\cdots p_{t+1}$, then
\begin{align}\label{eq4a}\frac{2^{\omega (n)}}{\varphi (n)}\le
\frac{2^t}{\varphi (p_1\cdots p_t)}\frac{p_1\cdots
p_t}{n}.\end{align}
\end{lemma}

\begin{proof} Since $n\ge p_1\cdots p_t$, there exists an integer
$s\ge t$ such that $p_1\cdots p_s\le n<p_1\cdots p_{s+1}$. Let the
standard prime factorization of $n$ be
$$n=q_1^{\alpha_1}\cdots q_k^{\alpha_k},$$
where $q_1, \dots , q_k$ are distinct primes and $\alpha_1, \dots
, \alpha_k$ are positive integers. Then $p_1\cdots p_k\le
q_1\cdots q_k\le n<p_1\cdots p_{s+1}$. It follows that $\omega
(n)=k\le s$. We also have
\begin{align*}\varphi (n)&=n(1-\frac 1{q_1})\cdots (1-\frac 1{q_k})\\
&\ge n(1-\frac 1{p_1})\cdots (1-\frac 1{p_k})\\
&\ge n(1-\frac 1{p_1})\cdots (1-\frac 1{p_s}).
\end{align*}
Hence \begin{align}\label{e2}\frac{2^{\omega (n)}}{\varphi
(n)}&\le \frac{2^{s}}{n(1-\frac 1{p_1})\cdots (1-\frac 1{p_s})}\nonumber\\
&= \frac{2^{s}}{\varphi (p_1\cdots p_s)} \frac{p_1\cdots
p_s}{n}\nonumber\\
&\le \frac{2^{s}}{\varphi (p_1\cdots p_s)} .
\end{align}
For $t<i\le s$, we have $p_i\ge 3$ and
$$\frac{2}{\varphi (p_i)}=\frac 2{p_i-1}\le 1.$$
It follows  from \eqref{e2} that \eqref{eq4} holds. It is clear
that \eqref{eq4a} follows from \eqref{e2} with $s=t$.

This completes the proof of Lemma \ref{lem8}.\end{proof}

\section{Proof of Theorem \ref{thm1}}

By Lemma \ref{lem1}, we have
\begin{align*}\pi (a, b)&=\sum_{k=1\atop (k,a)=1}^{a-2}(\pi (S; a, kb)-\pi (kb; a, kb)) +O(1)\\
&=\sum_{k=1\atop (k,a)=1}^{a-2}\Big( \frac 1{\varphi (a)}\int_{2}^{S}\frac 1{\log t} dt +O( S \exp (-D \sqrt{\log S}))\\
&\ \ \ -\frac 1{\varphi (a)}\int_{2}^{kb}\frac 1{\log t} dt -O( kb \exp (-D \sqrt{\log {kb}})) \Big)\\
&=\frac {\varphi (a)-1}{\varphi (a)}\int_{2}^{S}\frac 1{\log t} dt
-\sum_{k=1\atop (k,a)=1}^{a-2} \frac 1{\varphi
(a)}\int_{2}^{kb}\frac 1{\log t} dt\\
&\ \ \  +O( S \exp (-D \sqrt{\log S})).\end{align*} Given $1\le
k\le a-2$. Since
$$\int_{2}^{kb}\frac 1{\log t} dt=\frac{kb}{\log (kb)}-\frac 2{\log 2}+\int_{2}^{kb}\frac 1{\log^2 t} dt,$$
$$\frac{kb}{\log (kb)}=\frac{kb}{\log S} \frac{\log S}{\log (kb)},\quad \frac{\log S}{\log (kb)}\to 1 \text{ as } b\to +\infty ,$$
\begin{align*}0<\int_{2}^{kb}\frac 1{\log^2 t} dt &=\int_{2}^{\sqrt b}\frac 1{\log^2 t} dt +\int_{\sqrt b}^{kb}\frac 1{\log^2 t} dt\\
&<\frac{\sqrt b}{\log^2 2} +\frac{kb}{\log^2 \sqrt
b}\\
&<\frac{\sqrt S}{\log^2 2} +\frac{S}{\log S}\frac{4\log S}{\log^2
b}\end{align*} and
$$\frac{4\log S}{\log^2
b}\to 0  \text{ as } b\to \infty ,$$
 it follows that
$$\int_{2}^{kb}\frac 1{\log t} dt=\frac{kb}{\log S} (1+o(1))  \text{ as } b\to \infty .$$
Hence
\begin{align*}&\sum_{k=1\atop (k,a)=1}^{a-2} \frac 1{\varphi (a)}\int_{2}^{kb}\frac 1{\log t} dt\\
&=(1+o(1)) \frac b{\varphi (a) \log S} \sum_{k=1\atop (k,a)=1}^{a-2} k\\
&=(1+o(1)) \frac b{\varphi (a) \log S} (\frac 12 a\varphi (a) -(a-1))\\
&=(1+o(1)) \frac 12 \frac {ab}{\log S} - (1+o(1))\frac {b(a-1)}{\varphi (a) \log S} \\
&=(1+o(1)) \frac 12 \frac{a}{a-1} \frac {S}{\log S} - (1+o(1))\frac {S}{\varphi (a) \log S} \\
&= (1+o(1)) \left( \frac{a}{2(a-1)}-\frac 1{\varphi (a)}\right) \frac {S}{\log S}.
\end{align*}
Thus, by the prime number theorem, we have
\begin{align*}\pi (a, b)&=\frac {\varphi (a)-1}{\varphi (a)}
\int_{2}^{S}\frac 1{\log t} dt
-(1+o(1)) \left( \frac{a}{2(a-1)}
-\frac 1{\varphi (a)}\right) \frac {S}{\log S}\\
&\ \ \ \ +O( S \exp (-D \sqrt{\log S}))\\
&=(1+o(1)) \left( \frac {\varphi (a)-1}{\varphi (a)}-\frac{a}{2(a-1)}
+\frac 1{\varphi (a)}\right) \frac {S}{\log S}\\
&=\left( \frac 12-\frac 1{2(a-1)} +o(1)\right) \frac {S}{\log
S}\\
&=\left( \frac 12-\frac 1{2(a-1)} +o(1)\right) \pi (S).
\end{align*}

This completes the proof of Theorem \ref{thm1}.

\section{Proof of Theorem \ref{thm2}}

Since $(a-1)(b-1)=S+1$ and $3\le a<b$, it follows that $(a-1)^2<
S<(b-1)^2$. So $a< \sqrt S+1<b$.  It is clear that
\begin{equation}\label{e1}\pi (a, b)\ge \pi (S)-\pi (S-c)-\# \{ p\in \mathcal{P} : S-c< p\le S, p\notin
T( a, b)\} . \end{equation}

Now we divide into the following cases:

{\bf Case 1:} $a>6\cdot 10^4$.  Then $$a>2\cdot 3\cdot 5\cdot
7\cdot 11\cdot 13=30030,\quad S>(a-1)^2\ge 36\cdot 10^8.$$ Thus,
\begin{align*}\frac{\log (8a)}{\log S}&<\frac{\log (8(\sqrt S+1))}{\log S}\\
&=\frac 12+\frac{\log 8}{\log S}+\frac{\log (1+1/\sqrt S)}{\log S}\\
&<0.594504.\end{align*} If $a\ge 2\cdot 3\cdot 5\cdot 7\cdot
11\cdot 13\cdot 17$, then by Lemmas \ref{lem8},
\begin{align*}\frac{2^{\omega (a)}}{\varphi (a)}\le
\frac{2^7}{\varphi (2\cdot 3\cdot 5\cdot 7\cdot 11\cdot 13\cdot
17)}<0.0013889.\end{align*} If $6\cdot 10^4<a<2\cdot 3\cdot 5\cdot
7\cdot 11\cdot 13\cdot 17$, then by Lemmas \ref{lem8},
\begin{align*}\frac{2^{\omega (a)}}{\varphi (a)}<
\frac{2^{6}}{\varphi (30030)}\cdot \frac{30030}{6\cdot
10^4}<0.00556112.\end{align*} In view of  Lemma \ref{lem9}, we
have
\begin{align*}&\# \{ p\in \mathcal{P} : S-c< p\le S, p\notin
T( a, b)  \}\\
&<  \frac 14 \left( \frac 18 + \frac 34 \frac{1}{6\cdot 10^4}+
0.00556112 \right) \cdot  \frac{1}{ 1-0.594504} \cdot \frac{S}{\log S}\\
&<0.080503 \frac{S}{\log S} .
\end{align*}
It follows from \eqref{e1} and Lemma \ref{PiEstimate8numberic}
that
\begin{align*}\pi (a, b)&\ge \pi (S)-\pi (S-c)-\# \{ p\in
\mathcal{P} : S-c< p\le S, p\notin
T( a, b) \}\\
&> (0.121243 -0.080503) \frac{S}{\log S}>0.04\frac{S}{\log S}.
\end{align*}

{\bf Case 2:} $201\le a\le 6\cdot 10^4$. Then $S>(a-1)^2\ge 40000$
and
$$\frac{\log (8a)}{\log
S} < \frac{\log (8a)}{\log (a-1)^2}. $$ By Lemma \ref{lem9}, we
have
\begin{align*}&\# \{ p\in \mathcal{P} : S-c< p\le S, p\notin
T( a, b)  \}\\
&< \frac 1{4\varphi (a)}\Big(\sum_{\substack{1\le n\le (a+6)/8\\
\gcd (a,n)=1}} 1\Big) \Big( 1-\frac{\log (8a)}{\log
(a-1)^2}\Big)^{-1}\cdot \frac{S}{\log S}\\
&<0.11272 \frac{S}{\log S}.
\end{align*}
By \eqref{e1} and Lemma \ref{PiEstimate8numberic}, we have
\begin{align*}\pi (a, b)&\ge \pi (S)-\pi (S-c)-\# \{ p\in \mathcal{P} : S-c< p\le S, p\notin
T( a, b)  \}\\
&> (0.121243-0.11272 ) \frac{S}{\log S}>0.008\frac{S}{\log S}.
\end{align*}

{\bf Case 3:} $21\le a\le 200$ and $S\ge 40000$. Then
$$\frac{\log (8a)}{\log
S}\le \frac{\log (8a)}{\log 40000}. $$
 It follows from Lemma \ref{lem9} that
\begin{align*}&\# \{ p\in \mathcal{P} : S-c< p\le S, p\notin
T( a, b)  \} \\
&< \frac 1{4\varphi (a)}\Big(\sum_{\substack{1\le n\le (a+6)/8\\
\gcd (a,n)=1}} 1\Big) \Big( 1-\frac{\log (8a)}{\log
40000}\Big)^{-1} \cdot \frac{S}{\log S}\\
&<0.116219 \frac{S}{\log S}.
\end{align*}
By \eqref{e1} and Lemma \ref{PiEstimate8numberic}, we have
 \begin{align*}\pi (a, b)&\ge \pi (S)-\pi (S-c)-\# \{ p\in \mathcal{P} : S-c< p\le S, p\notin
T( a, b)  \}\\
&> (0.121243-0.116219 ) \frac{S}{\log S}>0.005\frac{S}{\log S}.
\end{align*}

{\bf Case 4: } $21\le a\le 200$ and $S<40000$. Noting that
$S+1=(a-1)(b-1)$, we have $b\le 40000/(a-1)+1$.  By Lemma
\ref{lem1}, for $21\le a\le 200$ and $202\le b\le 40000/(a-1)+1$
with $\gcd (a,b)=1$    we have
\begin{align*}\pi (a, b)&\ge  \# \{ (x, y)\in \mathbb{Z}^2 :  0\le x \le \frac b2-1,
0\le y\le \frac a2-1, ax+by\in \mathcal{P}\}\\
&\ge  \# \{ (x, y)\in \mathbb{Z}^2  : 0\le x \le 79,
0\le y\le \frac a2-1, ax+by\in \mathcal{P}\}\\
&=\# \{ (x, y)\in \mathbb{Z}^2  : 0\le x \le 79,
0\le y\le \frac a2-1, ax+by\in \mathcal{P}\}\\
&> 0.024  \frac{S}{\log S}.
\end{align*}
In view of Lemma \ref{lem1}, for $21\le a\le 200$ and $a< b\le
201$ with $\gcd (a,b)=1$ we have
\begin{align*}\pi (a, b)&\ge  \# \{ (x, y)\in \mathbb{Z}^2  :  0\le x \le \frac b2-1,
0\le y\le \frac a2-1, ax+by\in \mathcal{P}\}\\
&=\# \{ (x, y)\in \mathbb{Z}^2  :  0\le x \le \frac b2-1,
0\le y\le \frac a2-1, ax+by\in \mathcal{P} \} \\
&> 0.219 \frac{S}{\log S}.
\end{align*}

{\bf Case 5: } $3\le a\le 20$ and $b\ge 50 a^2$. By Lemma
\ref{pi(x;m,l)estimatesmall}, we have
\begin{align*}\pi (a, b)&=\pi (S; a, b)-\pi (b; a, b)\\
&>\frac 1{\varphi (a)} \frac{S}{\log S}-\frac 1{\varphi (a)}\frac{b}{\log b} \left( 1+\frac{5}{2\log b} \right)\\
&= \frac 1{\varphi (a)} \frac{S}{\log S}\left( 1- \frac{b}{S}
\frac{\log S}{\log b}\left( 1+\frac{5}{2\log b} \right)\right) .
\end{align*}
Since
$$\frac{b}{S}=\frac{b}{ab-b-a}=\frac 1{a-1-a/b}\le \frac 1{a-1-\frac 1{50a}},$$
$$\frac{\log S}{\log b}<\frac{\log (ab)}{\log b}=1+\frac{\log a}{\log b}\le 1+\frac{\log a}{\log (50a^2)},\quad
\frac 1{\log b}\le \frac 1{\log (50a^2)},$$ for $3\le a\le 20$ we
have
\begin{align*}\pi (a, b)&>\frac 1{\varphi (a)}
\left( 1-\frac 1{a-1-\frac 1{50a}}\left( 1+\frac{\log a}{\log (50a^2)}\right)
\left( 1+ \frac 5{2\log (50a^2)}\right) \right) \frac{S}{\log S}\\
&>0.05 \frac{S}{\log S}.\end{align*}

{\bf Case 6: } $3\le a\le 20$ and $a< b< 50 a^2$.
 By Lemma \ref{lem1}, for $3\le a\le 20$ and $a< b< 50a^2$ with
$\gcd (a, b)=1$ we have
\begin{align*}\pi (a, b)&\ge  \# \{ (x, y)\in \mathbb{Z}^2 :  x\ge 0,
y\ge 0, ax+by<S, ax+by\in \mathcal{P}\}\\
&\ge  \# \{ (x, y)\in \mathbb{Z}^2 : 0\le x \le 150,
y\ge 0, ax+by<S, ax+by\in \mathcal{P}\}\\
&=\# \{  (x, y)\in \mathbb{Z}^2  : 0\le x \le 150,
y\ge 0, ax+by<S, ax+by\in \mathcal{P}\} \\
&\ge 0.006 \frac{S}{\log S}.
\end{align*}

This completes the proof of Theorem \ref{thm2}.

\section*{Acknowledgments}

This work was supported by the National Natural Science Foundation
of China, Grant No. 12171243.

\end{document}